
\input amstex

\magnification=1200
\loadmsam
\loadmsbm
\loadeufm
\loadeusm
\UseAMSsymbols

\hsize=6.0truein
\hoffset=0.15truein
\vsize=9truein
\voffset=-0.2truein

\def\leftitem#1{\item{\hbox to\parindent{\enspace#1\hfill}}}

\def\boxit#1#2{\hbox{\vrule
	\vtop{%
	\vbox{\hrule\kern#1%
	\hbox{\kern#1#2\kern#1}}%
	\kern#1\hrule}%
	\vrule}}

\def\leaderfill{\leaders\hbox to 1em{\hss.\hss}\hfill}

\parskip=\medskipamount
\document

\input epsf

\centerline{\bf Cauchy Inequality and the Space of 
Measured Laminations, I}

\centerline{Feng Luo  and Richard Stong}

\bigskip

\S1. {\bf Introduction}

1.1.  This is the first 
of two papers addressing  a Cauchy type inequality for
the geometric intersection number between two 1-dimensional submanifolds
in a surface. 
As a consequence, we reestablish some of the basic results
in Thurston's theory of measured laminations. 
In this paper, we consider surfaces with non-empty boundary using
ideal triangulations.
In the sequel,  we establish the inequality
for closed surfaces using Dehn-Thurston coordinates.

1.2. Let us begin with a brief review  of Thurston's theory (see [Bo],
[FLP], [Mo], [PH], [Th1], [Th2] and others). Given a compact orientable surface 
$\Sigma$
with possibly non-empty boundary, a \it curve system \rm on $\Sigma$
is a proper 1-dimensional submanifold so that each component of it
is not null homotopic and not relatively homotopic into the boundary $\partial
\Sigma$. The space of all isotopy classes of curve systems on $\Sigma$ is
denoted by $CS(\Sigma)$. This space was introduced by Max Dehn in 1938 [De] who
called it the \it arithmetic field \rm of the topological surface. Given
two classes $\alpha$ and $\beta$ in $CS(\Sigma)$, their \it geometric intersection number
\rm $I(\alpha, \beta)$, is defined to be min\{ $|a \cap b|: a \in \alpha, b \in \beta$\}.
Thurston observed that the pairing I(,): $CS(\Sigma)$ $\times$ $CS(\Sigma)$ $\to \bold Z$
behaves like a non-degenerate ``bilinear'' form in the sense that (1) given
any $\alpha$ in $CS(\Sigma)$ there is $\beta$ in $CS(\Sigma)$ so that their intersection number
$I(\alpha, \beta)$ is non-zero, and (2) 
$I(k_1 \alpha_1, k_2 \alpha_2) = k_1k_2
I(\alpha_1, \alpha_2)$ for $k_i \in \bold Z_{\geq 0}$, $\alpha_i \in
CS(\Sigma)$ where $k_i \alpha_i$ is the collection of $k_i$ copies of
$\alpha_i$. 
In linear algebra,
given a non-degenerate quadratic form $\omega$ on a lattice
$L$ of rank $r$, one can form a completion of $(L , \omega)$ by
canonically embedding $L$ into $R^r$
so that the form $w$ extends continuously on $R^r$.
Thurston's construction is the exact analogy. 
Thurston's space of measured laminations on the surface $\Sigma$,
denoted by $ML(\Sigma)$ is defined to be the completion of the pair
$(CS(\Sigma), I)$ in the following sense. Given $\alpha$ in $CS(\Sigma)$, let
$\pi(\alpha)$ be the map sending $\beta$ to $I(\alpha, \beta)$. This gives an embedding
of $\pi: CS(\Sigma) \to \bold R^{CS(\Sigma)}$ where the target has
the the product topology.
The space $ML(\Sigma)$ is define to be the closure of 
$\bold Q_{>0} \times \pi(CS(\Sigma)) = \{ r \pi(x) : r \in \bold Q_{>0}, x \in $   $CS(\Sigma) \}$. 
Using the notion of train-tracks,
Thurston showed that $ML(\Sigma)$ is homeomorphic to a Euclidean space and 
that the intersection pairing I(,) extends to a continuous homogeneous
map from $ML(\Sigma)$$\times$$ML(\Sigma)$ to $\bold R$. See [Bo], [FLP], [Mo], [PH],
[Th2] and others for a proof of the first statement and [Bo],[Re] and
others for a proof of the continuity of the the extension. In fact
Bonahon generalized the intersection number to the geodesic currents and
proved the continuity of the pairing in the generalized situation.

One of our goals is to derive the continuity of the intersection pairing
on $ML(\Sigma)$ from a simple Cauchy inequality. This inequality also
implies that the space $ML(\Sigma)$ is a Euclidean space (see \S4).
Our main observation is that Thurston's theory of measured laminations
is based on a simple geometric fact that the intersection of
two distinct line segments in a disk consists of at most one point and
the intersection number being one is  detected by 
the relative locations of the end points of the line segments.

1.3. Given a compact surface $\Sigma$ with non-empty boundary,
recall that  \it an ideal triangulation \rm of a surface $\Sigma$ is a
maximal collection of pairwise disjoint, pairwise non-isotopic essential
arcs in $\Sigma$. Fix an ideal triangulation $t= t_1 \cup t_2 \cup ... \cup 
t_N$
of the surface, L. Mosher [Mo] proved that the space of closed
curve systems $CS_0(\Sigma)$ $=\{ [s] \in CS(\Sigma): \partial s =  \emptyset
\}$ can be parametrized using the ideal triangulation where
$[s]$ denotes the isotopy class of a submanifold $s$.
 Namely there is
an injective map $T:CS_0(\Sigma) \to \bold Z_{\geq 0}^N$
sending $\alpha$ to $(I(\alpha, [t_1]), ..., I(\alpha, [t_N])).$
Our main result is the following.

{\bf Theorem.} \it Suppose $t=t_1 \cup ... \cup t_N$ is an ideal
triangulation of a compact surface. Then for any three classes $\alpha,
\beta, \gamma \in CS_0(\Sigma)$, the following inequality holds:

$$|I(\alpha, \beta) - I(\alpha, \gamma)| \leq |\alpha||\beta - \gamma|$$
\noindent
where $|\alpha| = \sum_{i=1}^N I(\alpha, [t_i])$ and 
$|\beta -\gamma| = \sum_{i=1}^N |I(\beta, [t_i]) -I(\gamma, [t_i])|$.
Furthermore, the constant $1$ in the inequality is optimal. \rm

The basic idea of the proof is as follows. Given two classes
$\beta$, $\gamma$ so that $|\beta -\gamma| =n$, we produce a sequence of
1-dimensional submanifolds $\beta_i$ for $i=0,1, ..., n$ starting from $\beta$
and ending at $\gamma$ so that $|\beta_i - \beta_{i+1}| =1$. Thus,
we may assume that $|\beta - \gamma| =1$. We prove the theorem in
this special case by analyzing the surgery procedure 
relating $\beta$ and $\gamma$.

We remark that a similar result that
$|I(\alpha, \beta) - I(\alpha, \gamma)| \leq K |\alpha||\beta - \gamma|$
 was  obtained earlier by M. Rees [Re] using
train-tracks. The constant $K$ in her theorem is big and
depends on the train-tracks.

1.4. Let us derive the continuity of $I(,)$ on the space of compactly
supported measured laminations $ML_0(\Sigma)$ which is the
closure of $\bold Q_{\geq 0} \times \pi(CS_0(\Sigma))$ using the
above inequality.   One first extends the pairing $I( , )$ to 
$(\bold Q_{\geq 0} \times \pi(CS_0(\Sigma))^2$ by linearity
$I(k_1\alpha_1, k_2 \alpha_2) = k_1k_2 I(\alpha_1, \alpha_2)$. Thus
the inequality still holds for $\alpha$, $\beta$ and $\gamma$ in $\bold Q_{\geq 0}
\times CS_0(\Sigma)$. Since the product space $\bold R^{CS(\Sigma)}$
is metrizable, the continuity of the paring $I(,)$ on $ML_0(\Sigma)
\times ML_0(\Sigma)$
follows by showing that if $(\alpha_n, \beta_n) \in (\bold Q_{\geq 0}
\times CS_0(\Sigma))^2$ converges, then $I(\alpha_n, \beta_n)$
converges. Now since $\alpha_n$ and $\beta_n$ converge, both
lim$_n I(\alpha_n, [t_i])$ and lim$_n I(\beta_n, [t_i])$ exist for
all $t_i$. Thus, lim$_{n,m} |\alpha_n -\alpha_m| = 0$ and
lim$_{n,m} |\beta_n - \beta_m|$ = 0. By the inequality,
$|I(\alpha_n, \beta_n) - I(\alpha_m, \beta_m)|
\leq |I(\alpha_n, \beta_n) - I(\alpha_n, \beta_m)| +
|I(\alpha_n, \beta_m) -I(\alpha_m, \beta_m)|$
$\leq |\alpha_n| |\beta_n - \beta_m| + |\beta_m|| \alpha_n -\alpha_m|$
which conveges to 0.

As a consequence of the continuity, we see that the Cauchy inequality
still holds for $\alpha$, $\beta$ and $\gamma$ in $ML_0(\Sigma)$. 
Thus we deduce
a result of Mosher that each element $\alpha$ in $ML_0(\Sigma)$ is determined
by the $N$-tuple of intersection numbers 
 $T(\alpha) = (I(\alpha, [t_1])$,$ ...,$
$I(\alpha, [t_N]))$. Furthermore the Cauchy inequality implies
that the space $ML_0(\Sigma)$ is locally compact and
the map $T: ML_0(\Sigma) \to \bold R^N$ is proper. Indeed, if
a sequence $\alpha_n$ in $ML_0(\Sigma)$ is bounded under $T$, then
for any $\beta \in CS(\Sigma) $, the Cauchy inequality
implies that $I(\alpha_n, \beta) \leq  |T(\alpha_n)| |T(\beta)|$
is bounded in $n$ for each fixed $\beta$. Since there are at most
countably many $\beta$'s, by the standard diagonalization argument,
there is a subsequence $\alpha_{n_i}$ so that $I(\alpha_{n_i}, \beta)$
converges for all $\beta$. This simply says that $\{\alpha_n\}$
contains a  convergent subsequence. To see that $T$ is proper, 
we note that if $T(\alpha_n)$ converges to a point in $\bold R^N$,
then $T(\alpha_n)$ is bounded. Thus $\alpha_n$ contains a
convergent subsequence. This shows that $T$ is proper and 
$T: ML_0(\Sigma) \to \bold R^N$ is an embedding whose image is a closed
subset. 
In fact, Mosher showed that the image of $ML_0(\Sigma)$
is homeomorphic to a Euclidean space.  This result
will also be derived from theorem 1.3 in \S4.

1.5. In a forthcoming paper [LS], we shall establish an analogous
Cauchy inequality for the geodesic length of curve systems. Namely,
given a complete hyperbolic metric, let $l_d(\alpha)$ be
the sum of the lengths of the system of $d$-geodesics representing
$\alpha$. We prove that 
$|l_d(\alpha) - l_d(\beta)| \leq |d||\alpha -\beta|$ and 
$|l_d(\alpha) - l_{d'}(\alpha)| \leq |\alpha||d-d'|$
with respect to a fixed ideal triangulation of the surface. As a
consequence, we reestablish a result of Thurston that the
geodesic length function extends continuously to a function
on $Teich(\Sigma) \times ML_0(\Sigma) \to \bold R$. See [Bo]
for a written proof of this result.

1.6.  Part of the work is supported by the NSF.

\S2. {\bf Preliminaries}

2.1.
We begin by introducing some notations. Let $\Sigma$ = $\Sigma_{g,r}$ be a compact
orientable surface of genus $g$ with $r$ ($\geq 0$) many boundary
components. We shall assume that the Euler characteristic of
the surface $\Sigma$ is negative. Isotopies of the surface leave the
boundary invariant. Given a 1-submanifold $s$, we denote the isotopy
class of $s$ by $[s]$ and a small regular neighborhood of $s$ by
$N(s)$. The interior of a manifold $X$ will be denoted by $int(X)$.
The geometric intersection number $I([a], [b])$ will also be denoted
by $I(a, b)$ and $I([a], b)$. 

A simple loop (or a proper arc) in a surface $\Sigma$ is called \it trivial \rm
if it is null homotopic (or relatively homotopic into $\partial 
\Sigma$). A 1-submanifold in a surface $\Sigma$ is \it essential \rm
if each component of it is non-trivial. We shall enlarge the space of
 curve systems $CS(\Sigma)$ to the space of all isotopy classes of essential
1-submanifolds, denoted by $ES(\Sigma)$. The intersection pairing
$I(,)$ is defined similarly on $ES(\Sigma)$. Note that $CS(\Sigma)$ $\subset$
$ES(\Sigma)$ and $[\partial \Sigma] \in ES(\Sigma) -CS(\Sigma)$.

Our first result is  to give a parametrization of the space $ES(\Sigma)$ for
surfaces with non-empty boundary using ideal triangulations of the
surface. To achieve this, we begin by parametrizing  arc systems
on polygons.  

2.2. Let $P_n$ be an $n$-sided polygon. An \it arc \rm in $P_n$ is
a proper embedding of a closed interval into $P_n -$\{vertices of $P_n$\}. 
An arc in $P_n$ is called
\it trivial \rm if its end points either lie in one side of $P_n$ or
in two adjacent sides of $P_n$. An \it arc system \rm in $P_n$ is
a finite disjoint union of non-trivial arcs in $P_n$. Let
$ES(P_n)$  be the set of all isotopy classes of arc systems in $P_n$
where isotopies leave each side invariant. Given two classes $\alpha$
and $\beta$ in $ES(P_n)$, we define their intersection number to be
$I(\alpha, \beta) = min\{ |a \cap b| : a \in \alpha, b \in \beta\}$.
We say a non-trivial arc $s$ in $P_n$ is \it parallel \rm to a side
if one of the component of $P_n - s$ is a quadrilateral. 

We  first give a parametrization of $ES(P_6)$. Let the six
sides of the hexagon $P_6$ be $A_1, B_3, A_2, B_1, A_3, B_3$
labeled cyclically. A parametrization of $ES(P_6)$ using the $A$-sides
is as follows. Take $\alpha$ $=[a]$ in $ES(P_6)$. Let $x_i = I(\alpha, A_i) = 
|a \cap A_i|$ and $x_i'$ be the number of components of $a$ which are
parallel to $A_i$.  Evidently $x_i x_i'=0$. We call
$(x_1, x_2, x_3, x_1', x_2', x_3')$ the \it $t$-coordinate  \rm
of $\alpha$ with respect to the $A$-sides of the hexagon.

2.3. {\bf Lemma.} \it Let $\Delta =\{ (a_1, a_2, a_3)
\in \bold Z_{\geq 0}^3: a_i + a_j \geq a_k, $ for all $i\neq j \neq k 
 \neq i$\}.
The map $T: ES(P_6) \to $ \{$(x_1, x_2, x_3, x_1', x_2', x_3')
\in \bold Z_{\geq 0}^6: x_i x_i'=0$, if $(x_1, x_2, x_3) \in \Delta$
then $x_1 + x_2 + x_3$ is even\} sending an element
to its $t$-coordinate is a bijection. \rm

{\bf Proof.}  Clearly $T$ is well defined. To see that $T$ is onto, we construct
the arc system $a$ with a given vector $(x_1, x_2, x_3, x_1', x_2', x_3')$
as the coordinate according to the following five cases: (1) $ (x_1', x_2', x_3') =
(0,0,0)$, and $(x_1, x_2, x_3) \in \Delta$; (2) $(x_1', x_2', x_3')
=(0,0,0)$ and $(x_1, x_2, x_3) \notin \Delta$; (3) $x_i' =x_j' =0$,
$x_k' > 0$; (4) $x_i' = 0$ and $x_j'x_k' > 0$ and (5)
$x_1'x_2'x_3' > 0$. The corresponding arc systems are listed in
the figure 2.1  below.

 The arc system $a$ can be described as follows. Let
 $a_i$ (resp. $b_i$) be an arc parallel to $A_i$ (resp. $B_i$)
and $c_i$ be an arc joining $A_i$ to $B_i$. We use $kx$ to denote
$k$ parallel copies of an arc $x$. The in the case (1), $a = \cup_{k=1}^3
(\frac{x_i + x_j -x_k}{2}) b_k$; in the  case (2) say $x_k > x_i + x_j$,
then $a = x_i b_i \cup x_j b_j \cup (x_k - x_i -x_j) c_k$;
in the case (3), say $x_i \geq x_j$, then $a = x_k' a_k \cup
x_j b_j \cup (x_i-x_j) c_i$; in the case (4)  $a = x_j'a_j \cup
x_k'a_k \cup x_i c_i$ and in the last case (5) $a = \cup_{i=1}^3 x_i' a_i$. 
Since two non-trivial arcs are isotopic if and
only if their end points lands on the same set of sides on the polygon,
the map $T$ is injective. 
$\square$.

\midspace{0.1cm}
\centerline{\epsfbox{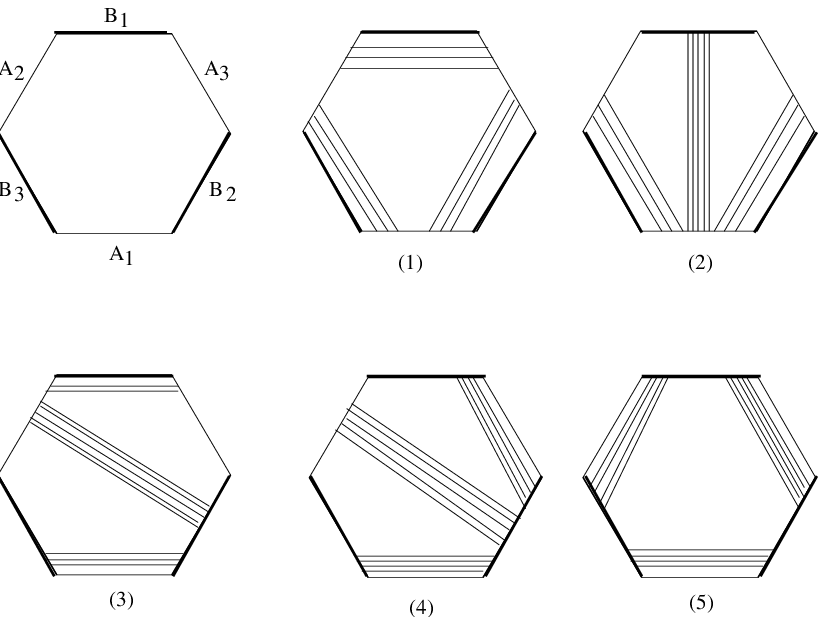}}
\centerline{ Figure 2.1}

2.4. Remark. The parametrization of the arc systems in
$P_6$ whose ends are in $A$-sides  is well known.

2.5. To parametrize the arc systems on any polygon $P_{2n}$ of even sides,
we use disjoint non-trivial arcs to decompose $P_{2n}$ into hexagons.
Let the $A$-sides of the hexagons correspond to the decomposing arcs.
Then a parametrization of $P_{2n}$ is given by taking the $t$-coordinates
of the hexagons with respect to the $A$-sides.  

2.6. One of the key ingredients in the proof of theorem 1.3 is to
understand the surgery procedure relating two elements in $ES(P_6)$
whose $t$-coordinates differ by a basis vector. For simplicity,
a class $\alpha$ in $ES(P_6)$ is called \it even \rm if all components
of its $t$-coordinates are even numbers.
We shall describe the surgery procedure relating two even arc systems
$\alpha$ and $\beta$ so that their $t$-coordinates 
$(x_1, x_2, x_3, x_1', x_2', x_3')$ and $(y_1, y_2, y_2, y_1', y_2', y_3')$ 
are related by $(x_1, x_2, x_3, x_1', x_2', x_3') =$  $(y_1, 
y_2, y_2, y_1', y_2', y_3') + (2,0,...,0)$.

Note that $x'_1=y'_1$. Since  $x_1>0$,  it follows that $x_1'=y_1'=0$.

Take a standard representative $a$ for $\alpha$. To obtain a
standard representative $b$ for  $\beta$, we perform the following
surgery operation on $a$. If $a$ contains arcs $b_2$ and $b_3$,
 we replace $a$ by $(a - b_2\cup b_3) \cup b_1$ to obtain $b$; if 
$a$ contains an arc parallel to $c_1$, then since $\alpha$ is even,
$a$ contains two copies of $c_1$. We replace $a$ by
$a  - 2 c_1$ to obtain $b$. In the remaining case, $a$ is disjoint
from either $c_2$ or $c_3$, say $a \cap c_2 = \emptyset$.
Since $x_1 \geq 2$, $a$ contains at least 2 copies of $b_3$.
In this case, replace $a$ by $(a - 2b_3) \cup 2c_2$ to obtain $b$.
Note that the arcs created lie in a small regular neighborhood of
the boundary and the arcs deleted.
See figure below for the illustration.

\midspace{0.1cm}
\centerline{\epsfbox{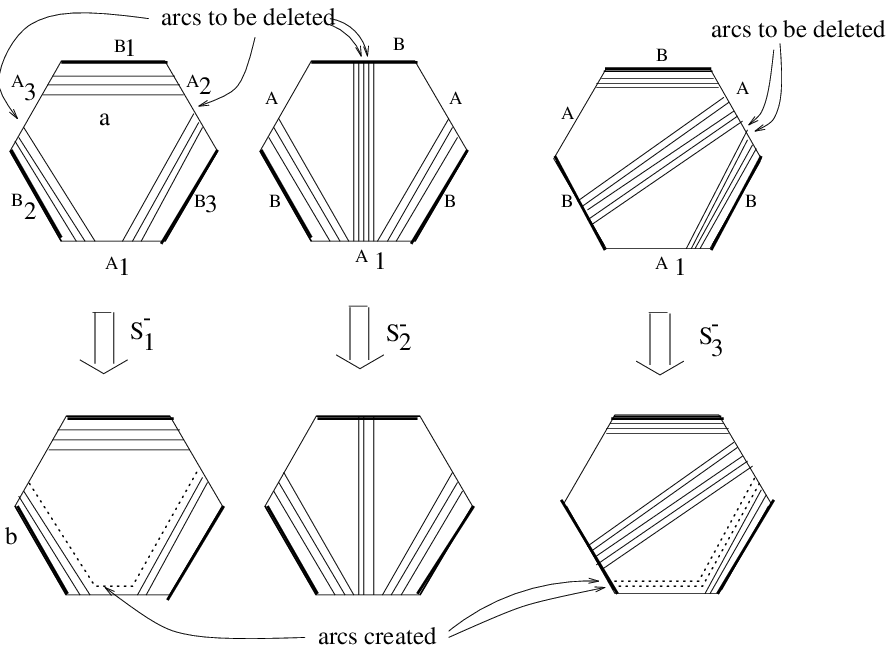}}
\centerline{Figure 2.2}
\midspace{0.1cm}
\centerline{\epsfbox{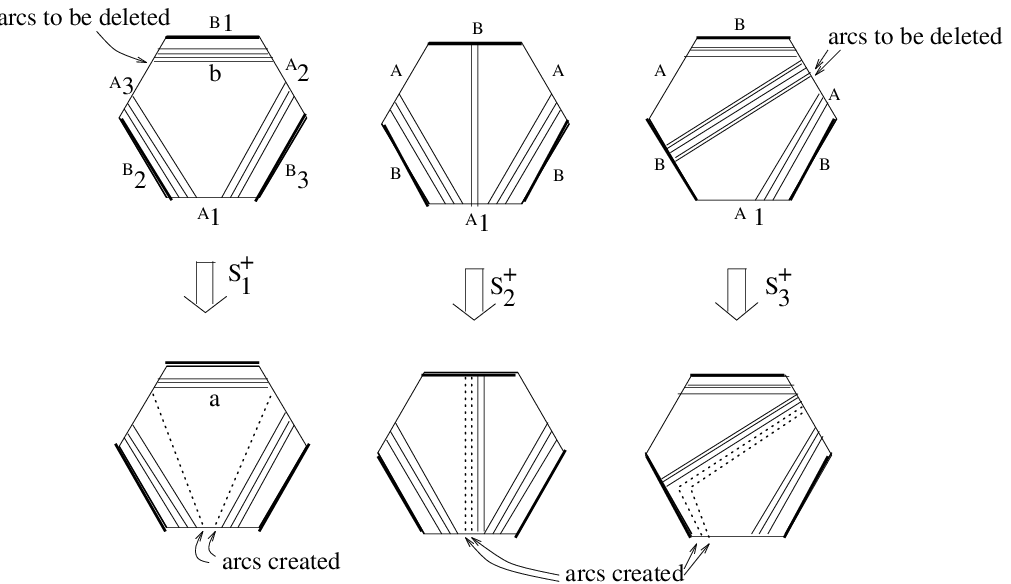}}
\centerline{Figure 2.3}

To obtain a standard representative of $a$ from $b$, we
perform the following surgery operation on $b$.
If $b$ contains some copies of $b_1$ but no $c_2$ or $c_3$,
replace $b$ by $(b -b_1)\cup b_2 \cup b_3$ to obtain $a$.
If $b$ contains no $b_1$, $c_2$ and $c_3$, replace $b$ by $b \cup 2c_1$
to obtain $a$. If $b$ contains some $c_2$ or $c_3$, say $c_2 \subset b$,
then $b$ contains even number of copies of $c_2$. Replace $b$ by
$(b -2c_2) \cup 2b_3$. 
 
\bigskip

\S3. {\bf Geometric Intersection Numbers on Surfaces with Boundary}

We  prove the main theorem 1.3 in this section. The basic idea of the
proof is as follows.  We  enlarge the classes of curve systems
to essential submanifolds and parametrize the space of all isotopy
classes of  essential
submanifolds by using an ideal triangulation. The reason for doing
so is due to the fact that given two isotopy classes of essential
submanifolds whose $t$-coordinates are distance-$k$ apart, there
is a sequence of $k+1$ essential 1-submanifolds starting and
ending at these two given classes so that the $t$-coordinates
of adjacent elements in the sequence are distance-1 apart. This
property does not hold for the space of curve systems. Thus, it
reduces the proof of the Cauchy inequality to the case of
$t$-coordinates being  distance-1 apart. By the surgery operation
relating distance-1 1-submanifolds, we prove the theorem.

3.1. We give a parametrization of the space $ES(\Sigma)$ as follows. Fix
a maximal collection $t= t_1 \cup ...\cup t_N$
of pairwise disjoint, non-isotopic essential arcs
(an ideal triangulation) of the surface $\Sigma$.  Thus the components
of $\Sigma - \cup_{i=1}^N int(N(t_i))$ are hexagons. Let the $A$-sides
of the hexagons correspond to $t_i$'s.
Given $\alpha$ in $ES(\Sigma)$, let  $t(\alpha)$ be the \it $t$-coordinate  \rm
of $\alpha$ which is the collection of $t$-coordinates of $\alpha$ in each hexagon. Namely, 
$t(\alpha) = (x_1, ..., x_N, x_1', ..., x_N')$ where $x_i = I(\alpha,
t_i)$ and $x_i'$ is the number of components of $\alpha$ equal to 
$[t_i]$. Clearly $x_ix_i'=0$.

3.2. {\bf Lemma.} (see also [Mo]) \it Fix an ideal triangulation $t$ of
$\Sigma$. Then the map $T: ES(\Sigma) \to X = \{(x_1, ..., x_N, x_1', ..., x_N')
\in \bold Z_{\geq 0}^N: x_ix_i'=0$, if $t_i, t_j$ and $t_k$ bound
a triangle and $(x_i, x_j, x_k) \in \Delta$, then $x_i + x_j + x_k$
is even\} sending an element to its
$t$-coordinate is a bijection.  In particular, the image of $T$ 
contains the set of even  vectors 
$ L = \{(x_1, ..., x_N, x_1', ..., x_N') \in (2 \bold Z_{\geq 0})^N :
 x_i x_i'=0\}$. \rm

Proof. To see that the map $T$ is onto,
take an element $(x_1, ..., x_N, x_1', ..., x_N')$ in the set $X$. Let
$H$ be a hexagonal component  of $\Sigma - \cup_{i=1}^N int(N(t_i))$ 
with three $A$-sides parallel to $t_i$, $t_j$ and $t_k$ 
(it may occur that $t_i = t_j$). By lemma 2.2, we construct an arc system in $H$ with 
the $t$-coordinate $(x_i, x_j, x_k, x_i', x_j', x_k')$. Now glue these
arc systems across $N(t_i)= t_i \times [-1, 1]$ by adding parallel
arcs $\{p_1, ..., p_n\} \times [-1, 1]$. We obtain a 1-submanifold $s$
properly embedded in $\Sigma$. By the construction, there are no
Whitney discs in $s \cup t$ and $s \cup \partial \Sigma$.
Thus the submanifold $s$ is essential and its $t$-coordinate is
the given vector  $(x_1, ..., x_N, x_1', ..., x_N')$. We call $s$
a \it standard representative\rm.  To see that the map $T$ is injective,
given $\alpha$ in $ES(\Sigma)$, choose a representative $a \in \alpha$ so that
$I(\alpha, t) = |a \cap t|$. Thus $a \cap H$ is an arc system in  
each hexagonal component of $\Sigma - \cup int(N(t_i))$. Since each non-trivial
arc in the quadrilateral $N(t_i)$  is parallel to a side, it follows
that $a$ is isotopic to a standard representative. It follows that
the map $T$ is injective.
$\square$

Having introduced the coordinate, let us now restate theorem 1.3 in
its most general form which we will prove.

3.3.
{\bf Theorem.} \it Suppose $t=t_1 \cup ... \cup t_N$ is an ideal
triangulation of a compact surface. Then for any three classes $\alpha,
\beta, \gamma \in ES(\Sigma)$ with $t$-coordinates
$(x_1, ..., x_N, x_1', ..., x_N')$,  
$(y_1, ..., y_N, y_1', ..., y_N')$ and
$(z_1, ..., z_N, z_1', ..., z_N')$, the following inequality holds:

$$|I(\alpha, \beta) - I(\alpha, \gamma)| \leq 2 |\alpha||\beta - \gamma|$$

\noindent
where $|\alpha| = \sum_{i=1}^N x_i + x_i' $ and 
$|\beta -\gamma| = \sum_{i=1}^N (|y_i-z_i| + |y_i' -z_i'|)$.
Furthermore, if
$\alpha$ is in $CS_0(\Sigma)$, then 
$$|I(\alpha, \beta) - I(\alpha, \gamma)| \leq |\alpha||\beta - \gamma|.$$

These inequalities are optimal. \rm

3.4. We now begin the proof of theorem 3.3 for classes in $ES(\Sigma)$. Since
the intersection pairing $I(,)$ is homogeneous, it suffices to prove
the inequality for $2 \alpha$, $2\beta$ and $2\gamma$ in $ES(\Sigma)$. The
$t$-coordinate of $2\alpha$ is an even vector in $L$. For simplicity,
we call a class $\alpha \in ES(\Sigma)$ \it even \rm if $T(\alpha) \in
L$ (see lemma 3.2). Thus it suffices to prove theorem 3.3 for even classes. 

Given two even vectors $u=(u_1, ...,$$u_{2N})$ and
$v=(v_1, ..., $$v_{2N})$ in $L$ so that their distance $|u-v|
=\sum_{i=1}^{2N} |u_i-v_i|$ is $2n$, there is a sequence of $n+1$
even vectors $w_j$, $j=0,..., n$ so that $w_0 =u$, $w_n=v$ and
$|w_{i+1} -w_i| =2$. 
Thus given two even classes $\beta$, $\gamma$ in $ES(\Sigma)$
so that $|\beta -\gamma| =2n$, by lemma 3.2, there exists a
sequence of $n+1$ even classes starting from $\beta$ and ending at $\gamma$
so that the adjacent elements are of distance-2 apart. Thus 
it suffices
to prove theorem 3.3 for even classes $\beta$ and $\gamma$ so
 that $|\beta -\gamma|=2$. Without loss of generality, we  
may assume that $T(\gamma) = T(\beta) +  (0,...,0,2,0,...,0)$,
i.e., $(z_1, ..., z_N, z_1', ..., z_N') = (y_1, ..., y_N, y_1', ..., y_N')
+ (0,...,0,2,0,...,0)$. We need to consider two cases: (1) $z_i'
=y_i'
+2$ and (2) $z_i = y_i + 2$ for some $i$.

In the first case that $z_i' = y_i'+2$, the class $\gamma$ is obtained
from $\beta$ by adding two copies of $[t_i]$. Thus $I(\alpha, \beta)
= I(\alpha, \gamma) -2 x_i$. The  inequality follows.

In the second case, 
we shall prove the Cauchy inequality by showing the following three
inequalities:

$$ I(\alpha, \beta) \leq I(\alpha, \gamma) + 4|\alpha| \tag 1$$ 
and
$$ I(\alpha, \gamma) \leq I(\alpha, \beta) + 2 |\alpha|. \tag 2$$

Furthermore, if $\alpha$ is a closed curve system in $CS_0(\Sigma)$,
we shall prove
$$I(\alpha, \beta) \leq I(\alpha, \gamma) + 2|\alpha|. \tag 3$$ 

Let us assume for simplicity that $i=1$. Let
$H_1$ and $H_2$ be the closures of the hexagonal 
components of $\Sigma - ( \cup_{i=2}^N
int(N(t)) \cup t_1)$ lying
on two sides of $t_1$ (it may occur that $H_1 = H_2 $). See figure 3.1.
If $H_1 \neq H_2$, then $H_1 \cap H_2 = t_1$
and $H_1 \cup H_2$ is an octagon. In this case we assume that $H_1 \cup
H_2$ is a convex octagon.

\midspace{0.1cm}
\centerline{\epsfbox{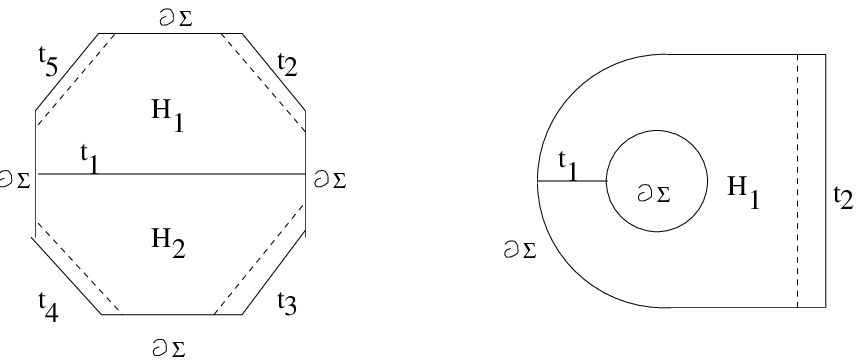}}
\centerline{Figure 3.1}

Take a standard representative $a$ of $\alpha$ so that $a \cap (H_1 \cup H_2)$
consists of straight line segments when $H_1 \cup H_2$ is an octagon.

3.5. To prove  $I(\alpha, \beta) \leq I(\alpha, \gamma) + 4|\alpha|$
$=I(\alpha, \gamma) + 2|\alpha||\beta - \gamma|$,
we find a standard representative $c$ of $\gamma$ so that
$|a \cap c| = I(\alpha, \gamma)$. By lemma 3.2, a standard
representative $b$ of $\beta$ can be constructed as follows. The submanifold
$b$ coincides with $c$ outside $H_1 \cup H_2 $.
In each of the hexagon $H_i - int(N(t_1))$, 
$b$ is obtained from $c$ by one of
the three  surgery operations $S_1^-, S_2^-$ or $S_3^-$
described in \S2.6, figure 2.2. Inside the regular
neighborhood $N(t_1)$ of $t_1$, the submanifold $b$ is obtained from
$c$ by a switch operation as shown in figure 3.2.

\midspace{0.1cm}
\centerline{\epsfbox{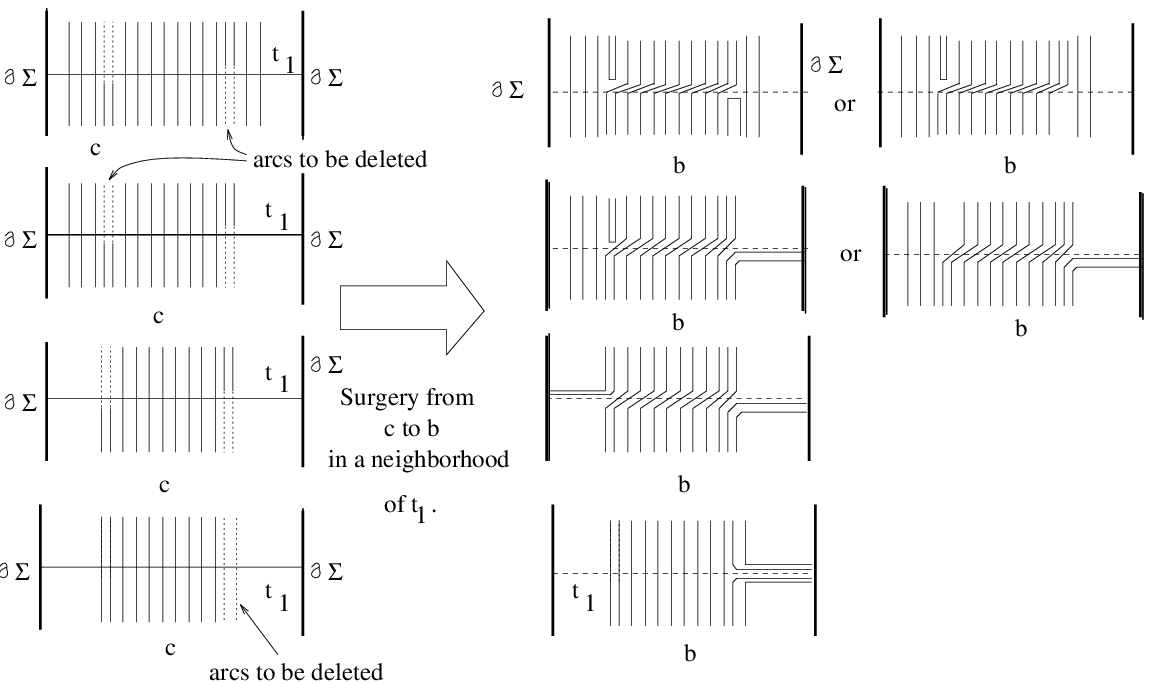}}
\centerline{Figure 3.2}

\midspace{0.1cm}
\centerline{\epsfbox{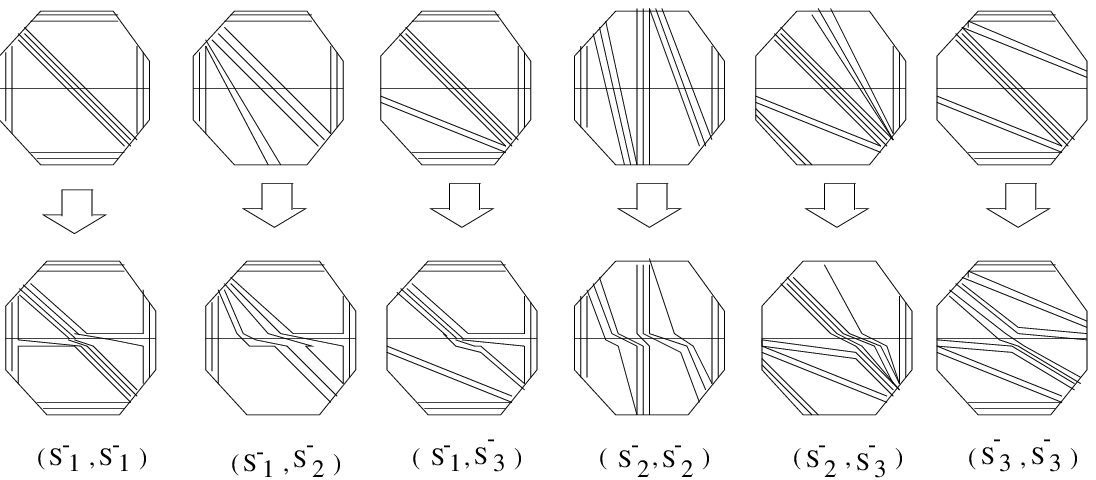}}
\centerline{All six surgeries relating the two submanifolds}
\centerline{Figure 3.3}

Note that the
dotted arcs in figure 3.2 indicate the arcs to be deleted in the surgery
or arcs whose end points are isotoped into $\partial \Sigma$. Also note
that the submanifold $b$ is obtained from $c$ by a surgery 
construction in $N(t_1)$ using at most four copies of $t_1$.
By definition, $I(\alpha, \beta) \leq |a \cap b|$. We estimate the
intersection number $|a \cap b|$ by considering the locations of
intersection points. 
Since the surgery operations in figure 2.2 show that the new arcs
created are in a small neighborhood of the deleted arcs and $t_1$,
$|a \cap b \cap (H_i -N(t_1))| \leq |a \cap c \cap (H_i -N(t_1))|$. 
By the construction, $|a \cap b \cap N(t_1)|
\leq 4 |a \cap t_1|.$
Thus $I(\alpha, \beta) \leq |a \cap b| \leq |a \cap c| + 4 |a|$
$= I(\alpha, \gamma) + 4 |\alpha|$.
If $\alpha \in CS_0(\Sigma)$, then the end points of $a \cap H_i$ are all
in the $A$-sides.  It follows that
$2|a \cap t_1| \leq |\alpha|$. Thus inequality (3) holds.

3.6. To prove $I(\alpha, \gamma) \leq I(\alpha, \beta) + 2 |\alpha|$
$= I(\alpha, \beta) +  |\alpha||\beta - \gamma|$,
we find a standard representative $b$ of $\beta$ so that $|a \cap b|
= I(\alpha, \beta)$ and
$b \cap (H_1 \cup H_2)$ consists of straight line segments  when
$H_1 \cup H_2$ is an octagon.
By lemma 3.2, 
a standard representative $c$ for $\gamma$ can be constructed as follows.
Let $c$ be the 1-submanifold coincide with $b$ outside 
$H_1 \cup H_2 $. In each of the hexagon $H_i -int(N(t_1)) $,
$c$ is obtained from $b$ by one of the three surgeries $S_1^+, S_2^+,
S_3^+$ described
in \S2.6 figure 2.3.
Inside the regular neighborhood $N(t_1)$ of $t_1$, the 1-submanifold
$c$ is obtained from $b$ by a switching operation as shown in figure  3.4.

\midspace{0.1cm}
\centerline{\epsfbox{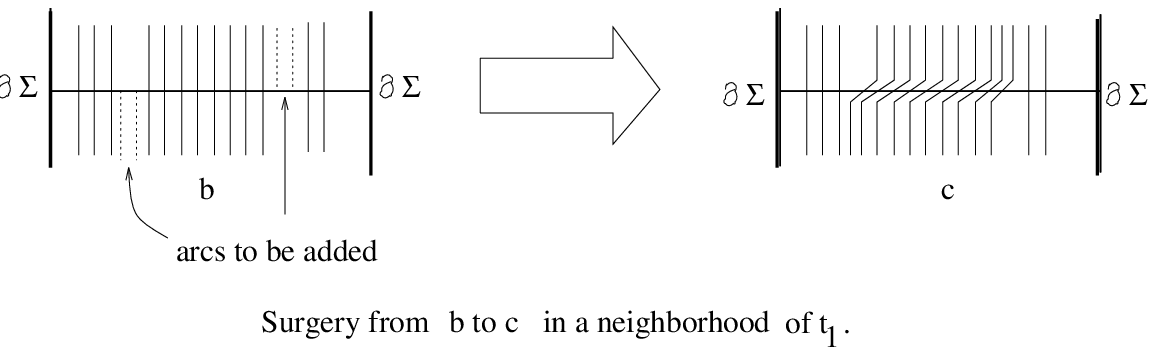}}
\midspace{0.1cm}
\centerline{Figure 3.4}

There are six types of surgeries relating $b$ to $c$. See figure 3.5.

\midspace{0.1cm}
\centerline{\epsfbox{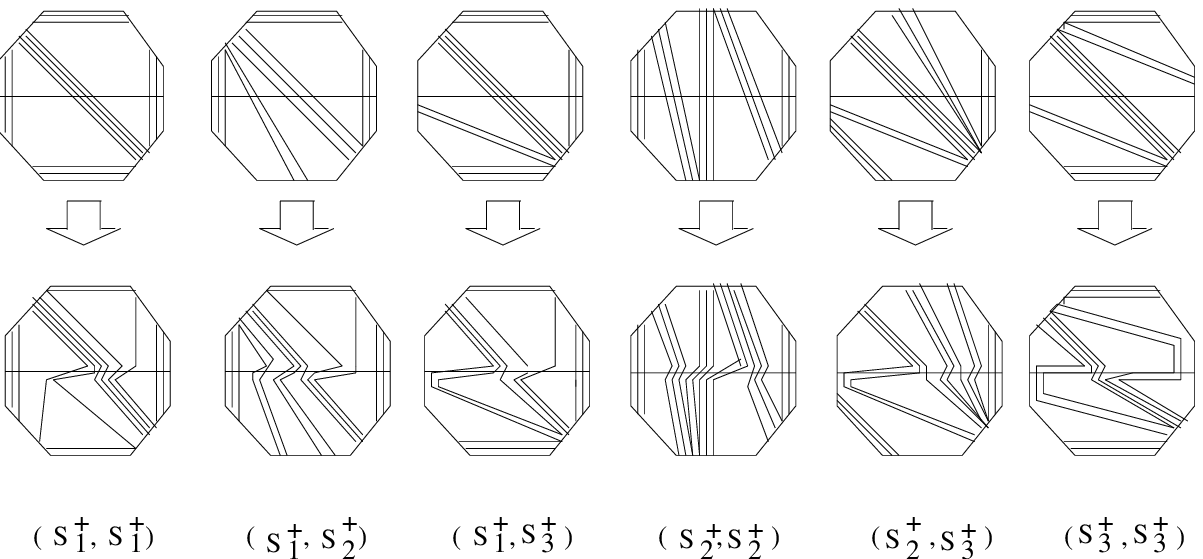}}
\midspace{0.1cm}
\centerline{Figure 3.5}

We need to consider the two cases that  $H_1 \neq H_2$ or $H_1 = H_2$
separately. Let us focus on the primary case that $H_1 \neq H_2$,
i.e., $H_1 \cup H_2$ is an octagon. In this case, let $c'$ be
the representative of $\gamma$ so that
1) $c'$ is equal to $c$ outside the octagon $H_1 \cup H_2$
and 2) $c' \cap (H_1 \cup H_2)$ consists of line segments. By
definition $I(\alpha, \gamma) \leq |a \cap c'|$. Thus to establish
inequality $(2)$, it suffices to show that inside the octagon
$H_1 \cup H_2$,

$$ |a \cap c' \cap (H_1 \cup H_2)| \leq |a \cap b \cap (H_1 \cup H_2)|
+ 2 |a|. \tag 4$$

To show $(4)$ and to save notation, 
consider only the local problem where $a$, $b$, and $c^{\prime}$ are
line segments in a single octagon $H_1 \cup H_2$. Thus we can write
merely $a$ instead of $ a \cap (H_1 \cup H_2)$, etc.
By the surgery construction relating $b$ and $c'$ (see figure 3.5),
we can express $b = b_1 \cup b_2$ and $c' = b_1 \cup b_2'$
as  disjoint unions of line segments so that
1) $b_2$ consists of parallel arcs crossing $t_1$,
2) $\partial b_2 \subset \partial b_2'$, $ |\partial b_2' -\partial b_2|
\leq 4$  and,
3) for each component $x$ of $b_2$ there is a component $x'$ of $b_2'$
so that $x$ and $x'$ share one end point and the other end points
of $x,x'$ bound an open interval containing at most one point of
$\partial b_2$.
See figure 3.6.

\midspace{0.1cm}
\centerline{\epsfbox{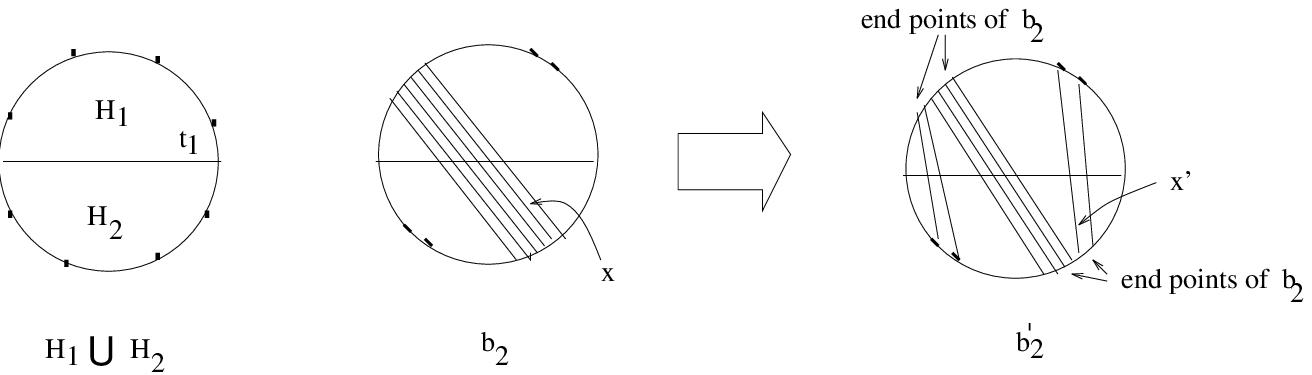}}
\midspace{0.1cm}
\centerline{Figure 3.6}

Now for each line segment $t$ in $H_1 \cup H_2$ whose end points
are in $\partial (H_1 \cup H_2)$, we have
$ |b_2' \cap t| - |b_2 \cap t| \leq 2$. Indeed, the intersection
of two line segments $s$ and $t$ in an octagon $H_1 \cup H_2$
is completely determined by the relative positions of their
end points $\partial s$ and $\partial t$. It follows that
inside $H_1 \cup H_2$,
$|a \cap c' \cap (H_1 \cup H_2) | \leq  |a \cap b \cap (H_1 \cup H_2)|
+ 2n$ where $n$ is the number of components of $a \cap (H_1 \cup H_2)$.
Since the number of components of an arc system having
$t$-coordinate $(x_1, x_2, x_3, y_1, y_2, y_3)$ inside
a hexagon is $\frac{x_1+x_2+x_3}{2} + y_1+y_2+y_3$, it follows
that $n \leq  |\alpha|$. (To see this there are several cases need to be
verified. Namely  one should discuss the cases where pair of $A$-sides
of $H_1 \cup H_2$ corresponds to a single $t_i$.) Thus inequality $(4)$
follows.

The second case that $H_1 = H_2$ is an annulus is simple. We simply
note that there are three surgeries relating $c$ to $b$ as shown in
figure 3.7. The three surgeries depend on the number $n$ of components
of arcs jointing $t_2$ to $\partial \Sigma$ in $H_1$. In the
first case, $n \geq 4$, in the second case $n=0$ and in the
last case $n=2$. In the first case $n \geq 4$, we remove four such
arcs and replace them by four arcs going around the boundary
component of $\partial \Sigma$. In the second case, we add two
parallel copies of the boundary components. In the last case
of $n=2$, we remove these two arcs and replace them by 
a parallel copy of the boundary component and an arc going around
the boundary.  Evidently inequality $(2)$ holds.

\midspace{0.1cm}
\centerline{\epsfbox{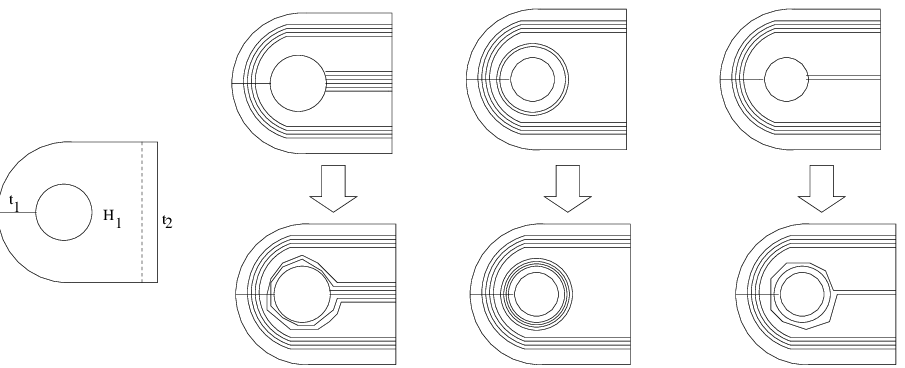}}
\midspace{0.1cm}
\centerline{Figure 3.7}

\S4. {\bf The Topology of the Space of Measured Laminations}

\bigskip

In this section, we derive the known fact  ([Mo], [Th]) that the 
spaces of all closed measured laminations $ML_0(\Sigma)$ is
homeomorphic to a  Euclidean space.

To see that $ML_0(\Sigma)$ is a Euclidean space, we fix
an ideal triangulation $t =t_1 \cup ... \cup t_N$ of the surface.
By theorem 1.3, we see that the map $T: ML_0(\Sigma) \to$
$ \bold R_{\geq 0}^N$ sending an element $\alpha$ to its
$t$-coordinate is an  embedding into a closed subset.
It remains to find the image of the map $T$. 
To this end, let us find the images under $T$ of the space
of all closed curve systems $CS_0(\Sigma)$.  Given
a $t$-coordinate $x =(x_1, ..., x_N)$ subject to the condition that 
when $t_i, t_j, $ and $t_j$ form the $A$-sides of a hexagonal
component of $\Sigma -\cup_r t_r$, then $(x_i, x_j, x_k) \in \Delta$,
one constructs an essential submanifold $s$ with  $x$ as its
$t$-coordinate by lemma 3.2. 
This essential submanifold $s$ is
a closed curve system if and only if  the submanifold $s$ contains
no loop parallel to $\partial \Sigma$. This is
the same as saying
that at least one of the hexagons
incident on $\partial_i\Sigma$ does not contain an arc parallel to the $B$-side
corresponding to $\partial_i\Sigma$, i.e., for each boundary
component $\partial_i \Sigma$,
$$\min_H \{ x_j+x_k-x_l \} = 0    \tag 8$$
where the minimum runs over all hexagons $H$ incident on $\partial_i\Sigma$ and
$H$ is formed by the arcs $t_j$, $t_k$, and $t_l$ with $t_l$ opposite to
a $B$-side in $\partial_i\Sigma$.
Suppose $r$ is the number of boundary components of the surface
$\Sigma$.
There are $r$ many equations $(8)$.
Thus we see that $CS_0(\Sigma)$ can be described as a finite union of regions,
each of which is described by integer coefficient linear equations (coming from (8)) in the $x_i$ and triangle inequalities saying that certain linear
combinations of the $x_i$ with integer coefficients are nonnegative. Thus the
set of rational solutions to these equations is dense in the set of real
solutions. This shows that the image 
$T(ML_0(\Sigma))$ is equal to the
subspace $S$ of $\bold R_{\geq 0}^N$ subject to $r$ equations $(8)$ and
the triangular inequalities: 

$$ x_j + x_k \geq x_l \tag 9$$
\noindent 
where {$t_j, t_k,$ $t_l$ form the $A$-sides of  a hexagonal component
of $\Sigma - \cup_{r=1}^N t_r$.

One may see the topological type of the  space defined by equations
$(8)$ and inequalities $(9)$ as follows. Let us make  a change of
variables by letting $y_i = ( x_j + x_k - x_l)/2$ in $(9)$.
Geometrically, $y_i$ is the number of copies of arcs parallel 
to the $B$-side of the hexagon corresponding to a boundary
component of $\partial \Sigma$. Then
$(9)$ becomes $(y_1, ..., y_M) \in \bold R_{\geq 0}^M$. The
equations $(8)$ become 
$$min \{ y_i | i \in B_j \} =0, \tag 10$$ 
where the index set $B_j$ consists of indices $i$ so that  $y_i$ is around the
$j$-th boundary component of the surface. 
Finally we  have a new set of equation
defined on each edge of the form $$y_i + y_j = y_k + y_l \tag 11 $$
for each edge $t_n$ of the ideal triangulation so that $y_i, y_j, y_k, y_l$
are adjacent to $t_n$.  These are exactly the switching equations in the
train-track dual to the ideal triangulation $t$ ([Mo], [Th]).
We claim that the equations $(10)$ and $(11)$ define a
space  $S$ in $\bold R^M$ homeomorphic to
a Euclidean space of dimension $6g-6+2r$. To this end,
consider the linear  subspace $V$ of $\bold R^M$ spanned by the vectors
$ \sum_{ i \in B_j} e_i$ 
where $e_i$ is the vector with $ y_i=1$ and all other $y_j=0$. Let $W$ be the
linear subspace defined by equations $(11)$.  
Let $P : \bold R^M \to R^M /V$ be the quotient map. We claim that
the restriction map $P|_S: S \to  P(W)$ is a homeomorphism.
Since $S$ is closed and locally compact, it suffices to show that
the restriction map $P|_S$ is one-to-one and onto. 
To see the map is onto,  given a vector 
$y =(y_1, ..., y_M )$ in $W$,  by adding  the vector 
$ - \sum_{j=1}^r  \sum_{i \in B_j}  {\text min}\{y_i | i \in B_j\} e_i $ 
to $y$, 
we see that the new vector is in the space $S$.  On the other hand, if
$y$ and $y'$ are two vectors in $S$ so that $y-y' \in V$, then by  looking
at the components around each boundary $\partial_j \Sigma$, we conclude that
$y=y'$.  This shows that the space $S$ and hence $ML_0(\Sigma)$ is homeomorphic
to a Euclidean space. To find the dimension of the Euclidean space, we note that
the two linear subspaces $W$ and $V$ intersect transversely  at $0$ in $R^M$.
Assuming this, since the dimension of $W$ is $6g+3r-6$ and the dimension of
$V$ is $r$, one finds the dimension of the quotient space to be $6g+2r -6$.

It remains to show that the subspaces $W$ and $V$ intersect transversely at $0$.
This
follows from a little bit of combinatorics, a linear combination of equations
of type  $(10)$ may be regarded as a linear combination of the duals to the
$t_i$ (suitably directed). Suppose a sum of these is a sum of 
equations of type  $(12)$. Then 
duals to consecutive $t_i$ around a boundary component must 
get weights which differ
by a constant.  But since the boundary edges cycle, this says that the duals
to the $t_i$ incident to a particular boundary component all get the same weight.
Since every boundary component is joined in a connected graph by the $t_i$,
we conclude that all duals get the same weight (up to sign for orientation). 
However looking at a single hexagon shows that the orientations cannot be
compatible unless all the weights are zero
and therefore all the weights are zero. Thus the only linear
combination which vanishes is the trivial one. This establishes the assertion.

Finally, we remark that the same argument shows that the closures of 
$\bold Q_{\geq 0} \times T(ES(\Sigma))$ and $\bold Q_{\geq 0} \times
T(CS(\Sigma))$ are Euclidean spaces. Their images in $\bold R^{2N}$
are given by $\{(x_1, ..., x_N, y_1, ..., y_N) \in \bold R^{2N}
|$  $x_i \geq 0$, $y_i \geq 0$ and $x_iy_i=0$ for all $i$\}
and $\{(x_1, ..., x_N, y_1, ..., y_N) \in \bold R^{2N}
|$  $x_i \geq 0$, $y_i \geq 0$, $x_iy_i=0$ for all $i$,
and equation $(8)$\}.

\centerline{\bf Reference}

[Bo] Bonahon, F.: Bouts des vari\'et\'es hyperboliques de dimension $3$. 
 Ann. of Math. (2) 124 (1986), no. 1, 71--158.

[De] Dehn, M.: Papers on group theory and topology. J. Stillwell (eds.).
 Springer-Verlag, Berlin-New York, 1987.

[FLP] Fathi, A., Laudenbach, F., Poenaru, V.: Travaux de Thurston sur les
surfaces. Ast\'erisque  66-67, Soci\'et\'e Math\'ematique de France, 1979.

[LS] Luo, F., Stong, R.: Cauchy inequality and geodesic length functions,
in preparation.

[Mo] Mosher, L.: Tiling the projective foliation space of a punctured surface.
Trans. Amer. Math. Soc. 306 (1988), 1-70.

[Pa] Papadopoulos, A.: Geometric intersection functions and
Hamiltonian flows on the space of measured foliations on a surface. Pacific
J. Math. 124 (1986), no. 2, 375--402.

[PH] Penner, R.,  Harer, J.: Combinatorics of train tracks.
Annals of Mathematics Studies, 125. Princeton University Press, Princeton,
NJ, 1992.

[Re] Rees, M.: An alternative approach to the ergodic theory of
measured foliations on surfaces. Ergodic Theory Dynamical Systems 1 (1981),
no. 4, 461-488. 

[Th1] Thurston, W.: On the geometry and dynamics of diffeomorphisms of
surfaces. Bul. Amer. Math. Soc. 19 (1988) no 2, 417-438.

[Th2] Thurston, W.: Geometry and topology of 3-manifolds, Princeton University
lecture notes, 1979.

\newpage

{\it address}

Feng Luo

Department of Mathematics

Rutgers University

New Brunswick, NJ 08854

fluo\@math.rutgers.edu

\bigskip

Richard Stong

Department of Mathematics

Rice University

Houston, TX  77005

stong\@math.rice.edu

\end